\documentclass{amsart}
\usepackage{amssymb,amsmath}

\newtheorem{theorem}{Theorem}

\newtheorem{corollary}{Corollary}

\newcommand{\RR}{{\mathbb R}}

\newcommand{\e}{\varepsilon}

\newcommand{\de}{\delta}
\newcommand{\del}{\partial}
\newcommand{\om}{\omega}
\newcommand{\Om}{\Omega}

\newcommand{{\loc}}{{\ell\mathrm oc}}

\newcommand{\calA}{{\mathcal A}}
\newcommand{\calB}{{\mathcal B}}
\newcommand{\calC}{{\mathcal C}}

\def\meanint{{\diagup\hskip -.42cm\int}}

\begin{document}

\title[C2-Differentiability of Solutions]
{Second-order Differentiability for Solutions of Elliptic Equations  in the Plane}
\author{Vladimir Maz'ya}
\address{Link\"oping University}
\email{vladimir.mazya@liu.se}
\author{ Robert McOwen}
\address{Northeastern University}
\email{mcowen@neu.edu}
\date{March 10, 2013}
\maketitle

\begin{abstract}
For a second-order elliptic equation of nondivergence form in the plane, we investigate conditions on the coefficients which imply that all strong solutions have first-order derivatives that are  Lipschitz continuous or differentiable at a given point. We assume the coefficients have modulus of continuity satisfying the square-Dini condition, and obtain additional conditions associated with a dynamical system that is derived from the coefficients of the elliptic equation.  Our results extend those of previous authors who assume the modulus of continuity satisfies the Dini condition. 


\end{abstract}

\smallskip\noindent

Let us consider an elliptic equation in nondivergence form
\begin{equation}\label{nondivergenceform}
a(x,y)\,u_{xx}+b(x,y)\,u_{xy}+c(x,y)\,u_{yy}=0  \quad\hbox{in}\  \Omega,
\end{equation}
where $\Omega$ is an open subset of ${\RR}^2$. 
Suppose that $u\in W^{2,2}(\Omega)$ is a strong solution of (\ref{nondivergenceform}); we want to know how regular
$u$ is in $\Omega$. This, of course, depends upon the smoothness of the coefficient functions $a$, $b$, and $c$. 
If we only assume the coefficients are bounded, then $u$ has first-order derivatives that are H\"older continuous in $\Omega$, i.e.\ $u\in C^{1,\alpha}(\Omega)$ where  $\alpha\in (0,1)$ depends on the coefficient bounds and the ellipticity constant (cf.\ \cite{N}); if the coefficients are continuous in $\Omega$, then $u\in C^{1,\alpha}(\Omega)$ for all $\alpha\in (0,1)$ (cf.\ \cite{ADN}).
On the other hand, if the coefficients are H\"older continuous in $\Omega$, or more generally
 if the coefficients are Dini-continuous in $\Omega$,  then it is well-known that $u\in C^2(\Omega)$ (cf.\ \cite{T}). In this paper, we want to find conditions on the coefficients, weaker than Dini continuity, under which $u$ will be second-order differentiable.

Let us assume that $\Omega$ contains the origin ${\bf 0}=(0,0)$ and focus on the differentiability at ${\bf 0}$.
Using ellipticity and a change of independent variables, we may arrange that $a({\bf 0})=1=c({\bf 0})$ and $b({\bf 0})=0$. 
Consequently, we  assume that the coefficients $a,b,c$ satisfy
\begin{equation}\label{nondivergence-coefficientconditions}
\sup_{|{\bf x}|=r}\left(|a({\bf x})-1|+|b({\bf x})|+|c({\bf x})-1|\right) \leq \om(r) \quad\hbox{as $r\to 0$, }
\end{equation}
where ${\bf x}=(x,y)$ and the modulus of continuity $\om$ is a continuous, nondecreasing function
for $0\leq r<1$ satisfying $\om(0)=0$. The coefficients being Dini continuous means that  (\ref{nondivergence-coefficientconditions}) holds with $\om(r)$ satisfying the Dini condition $\int_0^1 \om(r)\,r^{-1}\,dr<\infty$;
 H\"older continuity, of course, corresponds to the special case $\om(r)=C\,r^{\alpha}$ where $\alpha\in (0,1)$ and $C$ is a positive constant.  But we shall assume $\om(r)$ satisfies the more general square-Dini condition:
\begin{equation}
\int_0^1 \om^2(r)\,\frac{dr }{r}< \infty.
\label{eq:Sq-Dini}
\end{equation}

Given a solution $u\in W^{2,2}(\Omega)$ of  (\ref{nondivergenceform}), let us introduce the vector $U=(U_1,U_2)=(u_x,u_y)$. 
Using $u_{xx}=(U_1)_x$, $u_{xy}=(U_1)_y$, and $u_{yy}=(U_2)_y$, we can write (\ref{nondivergenceform}) as 
\[
a(x,y)\,(U_1)_x+b(x,y)\,(U_1)_y+(c(x,y)-1)(U_2)_y+(U_2)_y=0.
\]
If we differentiate this with respect to $x$ and use $(U_2)_{yx}=u_{yyx}=u_{xyy}=(U_1)_{yy}$ (where third-order derivatives are interpreted weakly), we obtain
\begin{equation}\label{U-equation}
(a(x,y)\,(U_1)_x)_x+(b(x,y)\,(U_1)_y)_x+((c(x,y)-1)(U_2)_y)_x+(U_1)_{yy}=0.
\end{equation}
Now we perform a similar calculation using $u_{xy}=(U_2)_x$ instead of $(U_1)_y$ and differentiating with respect to $y$ instead of $x$ to obtain 
\begin{equation}\label{V-equation}
(U_2)_{xx}+((a(x,y)-1)(U_1)_x)_y+(b(x,y)\,(U_2)_x)_y+(c(x,y)\,(U_2)_{y})_y=0.
\end{equation}
Putting (\ref{U-equation}) and (\ref{V-equation}) together as a second-order system, we obtain
\begin{equation}\label{second-order-system}
\left( \begin{pmatrix} a & 0 \\ 0 & 1 \end{pmatrix}  U_x \right)_x
+\left( \begin{pmatrix} 0 & 0 \\ (a-1) & b \end{pmatrix} U_x\right)_y
+\left(\begin{pmatrix} b & c-1 \\ 0 & 0 \end{pmatrix} U_y \right)_x
+\left(\begin{pmatrix} 1 & 0 \\ 0 & c \end{pmatrix} U_y\right)_y= {\bf 0}.
\end{equation}
 
 Now (\ref{second-order-system}) may look more complicated than (\ref{nondivergenceform}), but at least it is in divergence form:
\begin{equation}\label{divergenceformsystem}
(A_{11}U_x)_x+(A_{21}U_x)_y+(A_{12}U_y)_x+(A_{22}U_y)_y={\bf 0},
\end{equation}
where the $A_{ij}$ are $(2\times 2)$-matrices and $U\in W^{1,2}(\Omega,\RR^2)$ is a weak solution. Moreover, the matrices $A_{ij}$ are perturbations of $\delta_{ij} I$, where $I$ is the $2\times 2$ identity matrix, in that
\begin{equation}\label{eq:A_ij-(delta_ij)I}
\sup_{|{\bf x}|=r}|A_{ij}({\bf x})-\de_{ij} I|\leq \om(r) \quad\hbox{as $r\to 0$. }
\end{equation}
In this way, (\ref{divergenceformsystem}) is reminiscent of our work \cite{MM3} which considered the first-order differentiability of weak solutions to an elliptic equation in divergence form. Moreover, the first-order differentiability of $U$ corresponds to the second-order differentiability of $u$, so the conclusions of \cite{MM3} are just what we need here. However, the formulas of \cite{MM3} pertain to equations and not systems (whose coefficients are matrices so that products do not commute); consequently, they cannot be used directly in the present situation. Nevertheless, we can apply the methods of \cite{MM3} to (\ref{divergenceformsystem}).

The method of \cite{MM3} suggests that we find a first-order dynamical system on $0<t<\infty$ whose stability properties as $t\to\infty$ control the differentiability of the solutions of  (\ref{divergenceformsystem}). To derive the dynamical system, we first write
${\bf x}=r\,\theta$ where $r=|{\bf x}|$ and $\theta=(\theta_1,\theta_2)=(\cos\phi,\sin\phi)$ for $0\leq\phi<2\pi$. Then we write
\begin{equation}\label{def:U-decomp}
U(x,y)=U_0(r)+V_1(r)x+V_2(r)y+W(x,y),
\end{equation}
where $U_0$, $V_1$, and $V_2$ are given by mean integrals
\begin{equation}\label{def:U0,U1,U2}
\begin{aligned}
U_0(r)=&\,\meanint_{S^1} U(r\,\theta)\,d\phi,\\
V_1(r)=\meanint_{S^1} U(r\,\theta)\,\theta_1\,d\phi, & \quad
V_2(r)=\meanint_{S^1} U(r\,\theta)\,\theta_2\,d\phi,
\end{aligned}
\end{equation}
and $W(x,y)$ has zero spherical mean and first spherical moments:
\begin{equation}\label{W-conditions}
\meanint_{S^1} W(r\,\theta)\,d\phi=0=\meanint_{S^1} W(r\,\theta)\,\theta_1\,d\phi=\meanint_{S^1} W(r\,\theta)\,\theta_2\,d\phi.
\end{equation}
Similar to  \cite{MM3}, we shall show that the 4-vector function $\vec{V}(r)=(V_1(r),V_2(r))$ satisfies a dynamical  system  that depends on $W$, and $W$ satisfies a PDE that depends upon $\vec{V}$. Ultimately, we shall show that the behavior of $U_0$, $\vec{V}$, and $W$ are all controlled by the asymptotic behavior of solutions to the following first-order system:
 \begin{subequations}\label{DynSys}
\begin{equation}\label{DynSys:1}
\frac{d\phi}{dt}+R\,\phi=0 \quad\hbox{on}\ 0<t<\infty,
\end{equation}
where $r=e^{-t}$ and $R(e^{-t})$ is the $(4\times 4)$-matrix function defined on $0<t<\infty$ by
\begin{equation}\label{DynSys:2}
R(r):= \begin{pmatrix}
\overline{a}_1(r) & 0 & \overline{b}_1(r) & \overline{c}_1(r) \\
\overline{a}_2(r) & \overline{b}_2(r) & 0 & \overline{c}_2(r) \\
\overline{a}_2(r) & 0 & \overline{b}_2(r) & \overline{c}_2(r) \\
-\overline{a}_1(r) &  -\overline{b}_1(r)  & 0 & -\overline{c}_1(r)
\end{pmatrix},
\end{equation}
with coefficients  given by certain second spherical moments of the original coefficients:
\begin{equation}
\begin{aligned} 
\overline{a}_1(r):=\meanint_{S^1}a(r\theta)(\theta_2^2-\theta_1^2)\,d\phi, \qquad 
\overline{a}_2(r):=-2\,\meanint_{S^1}a(r\theta)\theta_1\theta_2\,d\phi, \\
\overline{b}_1(r):=\,\meanint_{S^1}b(r\theta)(\theta_2^2-\theta_1^2)\,d\phi, \qquad 
\overline{b}_2(r):=-2\,\meanint_{S^1}b(r\theta)\theta_1\theta_2\,d\phi, \\
\overline{c}_1(r):=\meanint_{S^1}c(r\theta)(\theta_2^2-\theta_1^2)\,d\phi, \qquad 
\overline{c}_2(r):=-2\,\meanint_{S^1}c(r\theta)\theta_1\theta_2\,d\phi.
\end{aligned}
\end{equation}
\end{subequations}

As in \cite{MM3}, the first-order regularity of solutions of (\ref{divergenceformsystem}) is determined by the stability properties of (\ref{DynSys}). In particular, we say that (\ref{DynSys}) is {\it uniformly stable} as $t\to\infty$ if for every $\e>0$ there exists a $\delta=\delta(\e)>0$ such that any solution $\phi$ of (\ref{DynSys:1}) satisfying $|\phi(t_1)|<\de$ for some $t_1>0$ satisfies $|\phi(t)|<\e$ for all $t\geq t_1$. We shall show that the following holds:

\begin{theorem}
If (\ref{DynSys}) is uniformly stable, then every weak solution $U\in W^{1,2}(\Omega,\RR^2)$  of (\ref{divergenceformsystem}) is Lipschitz continuous at ${\bf x}={\bf 0}$. 
\end{theorem}

Another important stability condition is that a solution of (\ref{DynSys}) be {\it asymptotically constant}, i.e.\ that $\phi(t)\to\phi_\infty$ as $t\to\infty$. As discussed in \cite{MM3}, this is actually independent of uniform stability, so we need to assume both conditions to conclude the differentiability of weak solutions. We shall show that the following holds:

\begin{theorem}
If (\ref{DynSys}) is uniformly stable and every solution is asymptotically constant, then every weak solution $U\in W^{1,2}(\Omega,\RR^2)$  of 
(\ref{divergenceformsystem}) is differentiable at ${\bf x}={\bf 0}$.
\end{theorem}

Recalling the derivation of (\ref{divergenceformsystem}) from (\ref{nondivergenceform}), these results yield the following:

\begin{theorem}
If (\ref{DynSys}) is uniformly stable, then every strong solution $u\in W^{2,2}(\Omega)$ of (\ref{nondivergenceform})
has first-order derivatives that are Lipschitz continuous at ${\bf 0}$. If, in addition, every solution of  (\ref{DynSys}) is asymptotically constant, then $u$ is second-order differentiable at ${\bf 0}$.
\end{theorem}

We can obtain analytic conditions on the matrix function $R$ that imply the desired asymptotic properties of (\ref{DynSys}). The simplest condition is 
\begin{equation}\label{RinL^1}
r^{-1}R(r)\in L^1(0,\e)\quad\hbox{for some $\e>0$,}
\end{equation}
 which guarantees that (\ref{DynSys}) is both uniformly stable and asymptotically constant; cf.\ \cite{C}. 

\begin{corollary}
If $R$ as in (\ref{DynSys:2}) satisfies (\ref{RinL^1}), then every strong solution $u\in W^{2,2}(\Omega)$ of (\ref{nondivergenceform}) is second-order differentiable at ${\bf 0}$.
\end{corollary}

\noindent
Analytic conditions weaker than (\ref{RinL^1}) can also be obtained. For example, if we introduce the symmetric matrix $S=-(R+R^t)/2$ and let $\mu(S)$ denote the largest eigenvalue of $S$, then it is shown in \cite{MM3} that
\begin{equation}\label{int-mu(S)}
\int_{r_1}^{r_2}\rho^{-1}\mu(S(\rho))\,d\rho<K \quad\hbox{for all $\e>r_2>r_1>0$}
\end{equation}
implies that (\ref{DynSys}) is uniformly stable; as a consequence, (\ref{int-mu(S)}) guarantees that every strong solution $u\in W^{2,2}(\Omega)$ of (\ref{nondivergenceform}) has first-order derivatives that are Lipschitz continuous at ${\bf x}={\bf 0}$. 
In addition, it is shown in \cite{MM3} that
\begin{equation}\label{R-int-RinL^1}
r^{-1}R(r)\int_0^r\rho^{-1}R(\rho)\,d\rho\in L^1(0,\e)
\end{equation}
implies that (\ref{DynSys}) is uniformly stable and asymptotically constant; as a consequence, (\ref{R-int-RinL^1}) guarantees that every strong solution $u\in W^{2,2}(\Omega)$ of (\ref{nondivergenceform}) is second-order differentiable at ${\bf x}={\bf 0}$. 

One may also consider  special cases to better understand the significance of the role of the dynamical system 
(\ref{DynSys}) in determining the regularity of strong solutions of (\ref{nondivergenceform}). In particular, let us assume that the coefficients $b$ and $c$ in (\ref{nondivergenceform}) satisfy
\begin{equation}\label{specialcase}
\overline{b}_i(r)=\overline{c}_i(r)=0\quad\hbox{ for $i=1,2$};
\end{equation}
this occurs, for example, when $b$ and $c$ are constant, or more generally if they depend only on $r$:
$b=b(r)$ and $c=c(r)$. In the case (\ref{specialcase}), we see that 
(\ref{DynSys:1}) decouples into three scalar equations
\[
\frac{d\phi}{dt}+\overline{a}_1\,\phi=0, \quad
\frac{d\phi}{dt}-\overline{a}_1\,\phi=0, \quad
\frac{d\phi}{dt}+\overline{a}_2\,\phi=0.
\]
But these are all of the form $\phi'+p(t)\phi=0$ which can be  solved using the integrating factor $\exp[\int^t p(\tau)d\tau]$. We conclude that the three scalar equations will be uniformly stable provided
\begin{equation}\label{examplecondn1}
\left|\int_s^t \overline{a}_1(\tau)\,d\tau\right|<K_1 \quad\hbox{and}\quad 
\int_s^t \overline{a}_2(\tau)\,d\tau >-K_2 \quad\hbox{for $t>s$ sufficiently large.}
\end{equation}
Moreover, the three scalar equations will be asymptotically constant provided
\begin{equation}\label{examplecondn2}
\begin{aligned}
\int_T^\infty \overline{a}_1(\tau)&\,d\tau \ \hbox{converges to a finite real number, and}\\
\int_T^\infty \overline{a}_2(\tau)\,d\tau  \  & \hbox{converges to an extended real number $>-\infty$.}
\end{aligned}
\end{equation}
Thus, when the coefficients $b,c$ satisfy (\ref{specialcase}) and the coefficient $a$ satisfies (\ref{examplecondn1}) and (\ref{examplecondn2}), then every strong solution of  (\ref{nondivergenceform}) will be second-order differentiable at ${\bf 0}$.

\section{Derivation of the Dynamical System}

A weak solution $U$ of (\ref{divergenceformsystem}) satisfies
\begin{equation}\label{divergenceformsystem-weak}
\int_\Om A_{ij}\,\del_iU \,\del_j\eta\,dxdy={\bf 0},
\end{equation}
for all $\eta\in C_0^\infty(\Omega)$, where we have used the Einstein summation convention of summing over repeated indices. To obtain the dynamical system (\ref{DynSys}), we begin by considering (\ref{divergenceformsystem-weak}) with different choices of test functions $\eta$. 

Taking $\eta$ to be a radial function $\eta(r)$ and using (\ref{def:U-decomp}), we obtain the following first-order ODE:
\begin{equation}\label{eq:ODE1}
\calA (r)\,U_0'+r\calB_1(r)\,V_1'+r\calB_2(r) \,V_2'+\Gamma_1(r)\, V_1+\Gamma_2(r) \,V_2+\Lambda[\nabla W](r)=0,
\end{equation}
where the $(2\times 2)$-matrices $\calA$, $\calB_k$, and $\Gamma_k$ are 
\begin{equation}\label{def:ODE1-matrices}
\begin{aligned}
\calA (r)&=\meanint_{S^1} A_{ij}(r\theta)\,\theta_i\theta_j\,d\phi =I_2+O(\om(r))\quad\hbox{as $r\to 0$} \\
\calB_k(r)&=\meanint_{S^1} A_{ij}(r\theta)\,\theta_k\theta_i\theta_j\,d\phi =O(\om(r)) \quad\hbox{as $r\to 0$} \\
\Gamma_k(r) & = \meanint_{S^1} A_{ik}(r\theta)\,\theta_i\,d\phi =O(\om(r)) \quad\hbox{as $r\to 0$},
\end{aligned}
\end{equation}
and the 2-vector $\Lambda$ is
\begin{equation}\label{def:ODE1-vector}
\Lambda[\nabla W](r)=\meanint_{S^1} (A_{i1}\theta_i W_x + A_{i2}\theta_i W_y)\, d\phi.
\end{equation}
Using Lemma 1 in \cite{MM3}, we can show that
\begin{equation}\label{est:ODE-vector}
|\Lambda[\nabla W](r)|\leq \om(r) \meanint_{S^1}|\nabla W|\,d\phi.
\end{equation}
Note that, although we are thinking of (\ref{eq:ODE1}) as an ODE, the coefficients are matrices and so it is really a system of two equations.

Taking $\eta=\eta(r)x$ and then $\eta=\eta(r)y$ in (\ref{divergenceformsystem-weak}), we obtain two second-order ODEs which we can put together as a second-order system:
\begin{equation}\label{eq:ODE2}
\begin{aligned}
-\left[ r^2\left( \begin{pmatrix} {\mathcal B}_1 \,U_0' \\ {\mathcal B}_2\, U_0' \end{pmatrix} 
+\begin{pmatrix}  {\mathcal A}_{11} & {\mathcal A}_{12} \\ {\mathcal A}_{21} & {\mathcal A}_{22} \end{pmatrix}
\begin{pmatrix} rV_1' \\ rV_2' \end{pmatrix}
+\begin{pmatrix}  {\calB}_{11} & {\calB}_{12} \\ {\calB}_{21} &{\calB}_{22} \end{pmatrix}
\begin{pmatrix} V_1 \\ V_2 \end{pmatrix}
+\begin{pmatrix} P_1[\nabla W] \\ P_2[\nabla W] \end{pmatrix}
\right) \right]' &\\
+\, r\left( \begin{pmatrix} \tilde \Gamma_1 U_0' \\ \tilde \Gamma_2 U_0' \end{pmatrix} + 
\begin{pmatrix}  \tilde{\mathcal B}_{11} & \tilde{\mathcal B}_{12} \\ \tilde{\mathcal B}_{21} & \tilde{\mathcal B}_{22} \end{pmatrix}
\begin{pmatrix} rV_1' \\ rV_2' \end{pmatrix}
+ \begin{pmatrix}  {\mathcal C}_{11} & {\mathcal C}_{12} \\ {\mathcal C}_{21} & {\mathcal C}_{22} \end{pmatrix}
\begin{pmatrix} V_1\\ V_2 \end{pmatrix}
+\begin{pmatrix} Q_1[\nabla W] \\ Q_2[\nabla W] \end{pmatrix}
\right)&=\begin{pmatrix}{\bf 0}\\{\bf 0}\end{pmatrix},
\end{aligned}
\end{equation}
where the ${\mathcal B}_j$ are defined above, but the other $(2\times 2)$-matrices are 
\begin{equation}\label{def:ODE2-matrices}
\begin{aligned}
\calA_{k\ell} (r)&=\meanint_{S^1} A_{ij}(r\theta)\,\theta_i\theta_j\theta_k\theta_\ell\,d\phi =\frac{1}{2}\,\delta_{k\ell}\,I_2+O(\om(r))\ \hbox{as $r\to 0$},\\
\calB_{k\ell} (r)&=\meanint_{S^1} A_{i\ell}(r\theta)\,\theta_i\theta_k\,d\phi =\frac{1}{2}\,\delta_{k\ell}\,I_2+O(\om(r))\ \hbox{as $r\to 0$},\\
\tilde{\calB}_{k\ell}(r)&=\meanint_{S^1} A_{ki}(r\theta)\,\theta_i\theta_\ell \,d\phi=\frac{1}{2}\,\delta_{k\ell}\,I_2+O(\om(r))
\ \hbox{as $r\to 0$},\\
\calC_{k\ell}(r)&=\meanint_{S^1} A_{k\ell}(r\theta)\,d\phi=\de_{k\ell}I_2+O(\om(r))\ \hbox{as $r\to 0$},\\
\tilde\Gamma_{k}(r) & =\meanint_{S^1} A_{ki}(r\theta)\,\theta_i\,d\phi=O(\om(r))\ \hbox{as $r\to 0$},
\end{aligned}
\end{equation}
and the $2$-vectors $P_1,\ P_2,\ Q_1,\ Q_2$ are given by
\begin{equation}\label{def:ODE2-vectors}
P_k[\nabla W](r)=\meanint_{S^1} A_{ij}(r\theta)\,\theta_i\theta_k \frac{\partial W}{\partial x_j}\,d\phi \quad \hbox{and} \quad 
Q_k[\nabla W](r)=\meanint_{S^1} A_{ki}(r\theta)\frac{\partial W}{\partial x_j}\,d\phi.
\end{equation}
As with (\ref{est:ODE-vector}), we can show
\begin{equation}\label{est:ODE2-vectorestimates}
|P_k[\nabla W](r)|, \ |Q_k[\nabla W](r)|\ \leq \om(r) \meanint_{S^1}|\nabla W|\,d\phi.
\end{equation}

We want to use (\ref{eq:ODE1}) to eliminate $U_0'$ from  (\ref{eq:ODE2}) and then identify the leading-order terms. Since ${\mathcal A}(r)=I_n+O(\om(r))$ as $r\to 0$, ${\mathcal A}(r)$ is invertible for small $r$ and we can write
\[
{\mathcal B}_i U_0'=-{\mathcal B}_i {\mathcal A}^{-1}\left({\mathcal B}_1rV_1'+{\mathcal B}_2rV_2'+\Gamma_1V_1+\Gamma_2V_2+\Lambda[\nabla W]\right).
\]
But  the coefficients of $rV_j'$ and $V_j$ in this expression are ``lower-order", i.e.\ 
\[
{\mathcal B}_i {\mathcal A}^{-1}{\mathcal B}_j, \ {\mathcal B}_i {\mathcal A}^{-1}{ \Gamma}_j = O(\omega^2(r))\quad
\hbox{as}\ r\to 0.
\]
So when we plug this into  (\ref{eq:ODE2}), it does not affect the leading order terms in $rV_j'$ and $V_j$. Similarly for replacing $\tilde \Gamma_i U_0'$ in  (\ref{eq:ODE2}).

Now let us make the substitution $r=e^{-t}$, so that $r\,d/dr=-d/dt$. Next, let us introduce
\[
\e(t)=\om(e^{-t}),
\]
and then write (\ref{eq:ODE2}) (after the elimination of $U_0$) as
\[
\left[ e^{-2t}\left( -{\bf A}\, \vec V_t  
+{\bf B}\, \vec V
+\vec P[\nabla W] +O(\e^2(t)) \right) \right]_t 
+\, e^{-2t}\left( 
-\widetilde{\bf B}\,
\vec V_t
+{\bf C}\,\vec V
+\vec Q[\nabla W]  + O(\e^2(t))
\right)=0,
\]
where $\vec V=(V_1,V_2)$, $\vec P=(P_1,P_2)$, and $\vec Q=(Q_1,Q_2)$ are 4-vectors and ${\bf A}$, ${\bf B}$, $\widetilde{\bf B}$, and ${\bf C}$ are the ($4\times 4$)-matrices
\[
{\bf A}=\begin{pmatrix}  {\mathcal A}_{11} & {\mathcal A}_{12} \\ {\mathcal A}_{21} & {\mathcal A}_{22} \end{pmatrix},
\quad
{\bf B}=\begin{pmatrix}  {\calB}_{11} & {\calB}_{12} \\ {\calB}_{21} &{\calB}_{22} \end{pmatrix},
\quad
\widetilde{\bf B}=\begin{pmatrix}  \tilde{\mathcal B}_{11} & \tilde{\mathcal B}_{12} \\ \tilde{\mathcal B}_{21} & \tilde{\mathcal B}_{22} \end{pmatrix},
\quad
{\bf C}= \begin{pmatrix}  {\mathcal C}_{11} & {\mathcal C}_{12} \\ {\mathcal C}_{21} & {\mathcal C}_{22} \end{pmatrix},
\]
and we have used $O(\e^2(t))$  to represent terms depending linearly on $\vec V_t$, $\vec V$, or $\Lambda[\nabla W]$, but with coefficients that are $O(\e^2(t))$ as $t\to\infty$. We can remove the factor $e^{-2t}$ to obtain
\begin{equation}\label{eq:ODE3}
\begin{aligned}
\left[  -{\bf A}\, \vec V_t  
+{\bf B}\, \vec V
+\vec P[\nabla W] 
+O(\e^2(t))  \right]_t  
&+ (2{\bf A}-\widetilde{\bf B}) \vec V_t \\
 +   ({\bf C}-2{\bf B}) \vec V
&-2\vec P[\nabla W]
+\vec Q[\nabla W]  + O(\e^2(t))
=0.
\end{aligned}
\end{equation}
However, this is still a second-order system, and we want to avoid differentiating the coefficient matrices, so let us convert it to a first order system by replacing the vector in the brackets in (\ref{eq:ODE3}) by a new 4-vector
\begin{equation}\label{def:U}
\vec U=-{\bf A}\, \vec V_t  +{\bf B}\, \vec V+\vec P[\nabla W] +O(\e^2(t)).
\end{equation}
We now have a first-order system in the 8-vector $(\vec V,\vec U)$:
\[
\begin{aligned}
\vec V_t-{\bf A}^{-1}{\bf B}\,\vec V+{\bf A}^{-1}\,\vec U&={\bf A}^{-1}\vec P[\nabla W] + O(\e^2) \\
\vec U_t+({\bf C}-\widetilde{\bf B}{\bf A}^{-1}{\bf B})\vec V+(\widetilde{\bf B}{\bf A}^{-1}-2{\bf I})\vec U&=
-\vec Q[\nabla W]+O(\e^2).
\end{aligned}
\]
where ${\bf I}$ is the ($4\times 4$) identity matrix. The coefficients of $\vec V$ and $\vec U$ behave as follows:
\[
-{\bf A}^{-1}{\bf B}\sim -{\bf I},
\quad {\bf A}^{-1}\sim 2{\bf I},
\quad {\bf C}-\widetilde{\bf B}{\bf A}^{-1}{\bf B}\sim \frac{1}{2}{\bf I},
\quad \widetilde{\bf B}{\bf A}^{-1}-2{\bf I}\sim -{\bf I}.
\]
where $\sim$ means differs by a term that is $O(\e)$ as $t\to\infty$.
Consequently, let us rewrite the first-order system as
\begin{equation}\label{ODE4}
\frac{d}{dt}\begin{pmatrix} \vec V \\ \vec U \end{pmatrix}
+ {\bf M}(t)\begin{pmatrix} \vec V \\ \vec U \end{pmatrix}
=\begin{pmatrix} {F_1(t,\nabla W)} \\ {F_2(t,\nabla W)} \end{pmatrix}
\end{equation}
where the ($8\times 8$)-matrix-valued function ${\bf M}(t)$ is of the form
\[
{\bf M}(t)={\bf M}_\infty + {\bf S}_1(t)+{\bf S}_2(t),
\]
with a constant matrix
\[
{\bf M}_\infty=\begin{pmatrix} -{\bf I} & 2\,{\bf I} \\ \frac{1}{2}\,{\bf I} & - {\bf I} \end{pmatrix}.
\]
The variable coefficient matrices ${\bf S}_1$ and ${\bf S}_2$ satisfy
\[
\begin{aligned}
 {\bf S}_1(t)=\begin{pmatrix} {\bf I}-{\bf A}^{-1}{\bf B} & {\bf A}^{-1}-2{\bf I}\\ {\bf C}-\widetilde{\bf B}{\bf A}^{-1}{\bf B}-\frac{1}{2}{\bf I}
& \widetilde{\bf B}{\bf A}^{-1}-{\bf I}\end{pmatrix}=O(\e(t))\quad \hbox{as}\ t\to\infty 
\end{aligned}
\]
and ${\bf S}_2=O(\e^2(t))$ as $t\to\infty$. The right-hand side of (\ref{ODE4}) satisfies
\[
|F_i(t,\nabla W)|\leq   \e(t) \meanint_{S^1}|\nabla W|\,d\phi.
\]

In order to analyze (\ref{ODE4}), as in \cite{MM3} we introduce a change of variables
\begin{equation}\label{def:phi,psi}
\begin{pmatrix} \vec V \\ \vec U \end{pmatrix} = {\bf J} \begin{pmatrix} \phi \\ \psi \end{pmatrix},
\end{equation}
where the matrix
\[
{\bf J}=\begin{pmatrix} 2\,{\bf I} & 2\,{\bf I} \\ {\bf I} & -{\bf I} \end{pmatrix}
\]
diagonalizes ${\bf M}_\infty$, i.e.\ ${\bf J}^{-1}{\bf M}_\infty {\bf J}=\hbox{diag}(0,0,0,0,-2,-2,-2,-2)$. We find that $(\phi,\psi)$ satisfies a dynamical system of the form
\begin{equation}\label{ODE5}
\frac{d}{dt} \begin{pmatrix} \phi \\ \psi \end{pmatrix} +
 \begin{pmatrix} 0 & 0 \\ 0 & -2\,{\bf I} \end{pmatrix}\begin{pmatrix} \phi \\ \psi \end{pmatrix}
 +{\bf R}(t) \begin{pmatrix} \phi \\ \psi \end{pmatrix} =G(t,\nabla W),
\end{equation}
where
\[
{\bf R}(t)=\begin{pmatrix} R_1(t) & R_2(t) \\ R_3(t) & R_4(t) \end{pmatrix},
\]
with
\[
R_1(t)\approx \frac{1}{4}{\bf A}^{-1}-\frac{1}{2}{\bf A}^{-1}{\bf B}+{\bf C}-\widetilde{\bf B}{\bf A}^{-1}{\bf B}+\frac{1}{2}\widetilde{\bf B}{\bf A}^{-1}-{\bf I},
\]
where $\approx$ means differs by a term that is $O(\e^2(t))$ as $t\to\infty$. The right-hand side of (\ref{ODE5}) satisfies
\begin{equation}\label{est:G}
|G(t,\nabla W)|\leq   \e(t) \, \meanint_{S^1}|\nabla W|\,d\phi.
\end{equation}
Estimates on $W$ that we shall discuss in the next section together with the stability theory presented in Section 2 of \cite{MM3} show that the stability of (\ref{ODE5}) is determined by that of 
\begin{equation}
\frac{d\phi}{dt}+R_1(t)\phi =0,
\end{equation}
so we need to determine the asymptotic behavior of $R_1$. But to do this, let us write
\[
{\bf A}= \frac{1}{2}({\bf I}+{\bf A}_0),  \quad
{\bf B}=\frac{1}{2}({\bf I}+{\bf B}_0), \quad \widetilde{\bf B}= \frac{1}{2}({\bf I}+\widetilde{\bf B}_0), 
\quad {\bf C}= {\bf I}+{\bf C}_0,
\]
where $|{\bf A}_0|,|{\bf B}_0|, |\widetilde{\bf B}_0|, |{\bf C}_0|=O(\e(t))$ as $t\to \infty$. Also note that
\[
\quad {\bf A}^{-1}\approx 2({\bf I}-{\bf A}_0).
\]
Using these we can simplify $R_1$ to obtain
\[
R_1\approx {\bf C}_0-{\bf B}_0={\bf C}-2{\bf B},
\]
and after a careful calculation we obtain the formula given in (\ref{DynSys:2}).

\section{Proofs of Theorems 1 and 2}

Since we are only interested in the behavior of our weak solution near ${\bf 0}$, 
we may assume $\Omega=B_\e({\bf 0})$ with $\e>0$ chosen small enough to make 
\begin{equation}\label{delta-conditions}
\int_0^\e r^{-1}\om(r)\,dr<\de \quad\hbox{and}\quad \om(\e)<\de,
\end{equation}
with $\de>0$ as small as we like. In fact, for any $p\in (1,\infty)$ we can choose $\delta=\delta(p)>0$ in
(\ref{delta-conditions}) small enough that the small oscillation condition on the coefficients (\ref{nondivergence-coefficientconditions}) ensures that $\nabla U\in L^p_{\loc}(\Omega)$; cf.\ Corollary 6.2 in \cite{MMS}.
Henceforth, we pick $p>2$ and choose $\epsilon$ small enough that $\nabla U\in L^p_{\loc}(\Omega)$.
But by rescaling the independent variables, we may arrange $\e>1$, so we may assume that our weak
solution $U$ of (\ref{divergenceformsystem}) satisfies
\begin{equation}\label{grad-U_is_Lp}
\nabla U\in L^p(\Omega)\quad\hbox{where}\ p>2 \ \hbox{and}\ \Omega=B_1(0).
\end{equation}
In particular, by Sobolev's inequality we know that $U$ is continuous in $\Omega$.

For our analysis, it is useful to consider (\ref{divergenceformsystem}) on all of $\RR^2$, so we extend the matrices $A_{ij}$ to all of $\RR^2$ by
\[
{ A}_{ij}=\de_{ij} I \quad\hbox{for $|x|>1$.}
\]
We also extend our modulus of continuity $\om$ to $(0,\infty)$ by $\om(1)$ for $r>1$.
It will also be useful to introduce the $L^p$-mean of a function over the annulus $A_r=\{x:r<|{\bf x}|<2r\}$:
\[
M_p(f,r)=\left(\meanint_{A_r}|f(x)|^p\,dx\right)^{1/p}.
\]
To control growth of the first derivatives of functions, we introduce
\[
M_{1,p}(f,r)=rM_p(\nabla f,r)+M_p(f,r).
\]

Let us introduce a smooth cut-off function $\chi(r)$ that is $1$ for $0\leq r\leq 1/4$ and $0$ for $r\geq 1/2$. We find that $\chi(r) U(x,y)$   satisfies 
\[
(A_{11}(\chi U)_x)_x+(A_{21}(\chi U)_x)_y+(A_{12}(\chi U)_y)_x+(A_{22}(\chi U)_y)_y=
F_0+(F_1)_x+(F_2)_y
\]
where
\[
F_0=A_{11}\chi'\theta_1 U_x+A_{21}\chi'\theta_2U_x+A_{12}\chi'\theta_1U_y+A_{22}\chi'\theta_2U_y,
\]
\[
F_1=\chi'(A_{11}\theta_1+A_{12}\theta_2)U
\quad\hbox{and}\quad 
F_2=\chi'(A_{21}\theta_1U+A_{22}\theta_2)U.
\]
Using (\ref{divergenceformsystem-weak}) with $\eta=\chi$, we see that 
\begin{equation}\label{int-F0=0}
\int_{\RR^2} F_0\,dxdy=0.
\end{equation}
Since we are interested in the behavior near $x=0=y$ where $U$ and $\chi U$ agree, we can simply assume that $U$ is supported in $r\leq 1/2$ and satisfies
\begin{equation}\label{nonhomogeneousdivergenceformsystem}
\del_i(A_{ij}\del_jU)=
F_0+\del_i(F_i),
\end{equation}
where $F_0,F_1,F_2\in L^p(\RR^2)$ are supported in $1/4\leq r  \leq 1/2$ and $F_0$ satisfies (\ref{int-F0=0}).
Of course, we now must replace (\ref{divergenceformsystem-weak}) by
\begin{equation}\label{nonhomogeneousdivergenceformsystem-weak}
\int_{\RR^2} A_{ij}\,\del_j U\,\del_i\eta \,dxdy=
\int_{\RR^2}(F_i\del_i\eta- F_0)\,dxdy,
\end{equation}
for all $\eta\in C_0^\infty(\Omega)$. 

At this point we observe that (\ref{nonhomogeneousdivergenceformsystem}) with (\ref{int-F0=0}) for the vector function $U$ is identical with (51ab) in \cite{MM3} for the scalar function $u$. This means that we can repeat the analysis of \cite{MM3} to connect the stability of the dynamical system (\ref{DynSys:2}) with the regularity of our weak solution. 
We do not want to repeat all of the details here, but let us give an outline of the argument. 

To begin with, we recall the decomposition $U=U_0+V_1 x+ V_2 y + W$ in (\ref{def:U-decomp}). We have shown that $\vec V$ satisfies a dynamical system (\ref{ODE4}) that depends on $\nabla W$, so we need to know $\nabla W$ is sufficiently well-behaved in order to obtain estimates for $\vec V$. This is done by showing that $W$ satisfies a PDE that depends on $\vec V$. To derive the PDE for $W$, we introduce
\begin{equation}
\Omega_{ij}=A_{ij}-\de_{ij} I,
\end{equation}
which satisfies
\[
|\Omega_{ij} | \leq \om(r) \quad\hbox{for}\ 0<r<1\quad\hbox{and}\quad \Omega_{ij}=0 \quad \hbox{for}\ r\geq 1.
\]
We also introduce for $f\in L^1_{\loc}(\RR^2\backslash\{0\})$ 
\begin{equation}
f(r\theta)^\perp=f(r\theta)-Pf(r\theta),
\end{equation}
where $P$ is the projection of $f$ onto the functions on $S^1$ spanned by $1,\theta_1,\theta_2$:
\[
\begin{aligned}
Pf(r\theta)&=c_0(r)+c_1(r)\theta_1+c_2(r)\theta_2, \ \hbox{where} \\
c_0(r)=\meanint_{S^1}&f(r\theta)\,d\phi\quad\hbox{and}\quad c_i(r)=2\,\meanint_{S^1} \theta_i\,f(r\theta)\,d\phi.
\end{aligned}
\]
Notice that 
\[
P[\Delta (U_0+V_1x+V_2y)]=\Delta (U_0+V_1x+V_2y)
\quad\hbox{and}\quad P[\Delta W]=0,
\]
so $W$ satisfies the following perturbation of Laplace's equation on $\RR^2$:
\begin{equation}\label{W-PDE}
\Delta W + \left[\del_i(\Omega_{ij}\del_j U_0)\right]^\perp + \left[\del_i(\Omega_{ij}\del_j(V_kx_k))\right]^\perp
+\left[\del_i(\Omega_{ij}\del_jW)\right]^\perp=[F_0+\del_i(F_i)]^\perp.
\end{equation}

Now we simultaneously consider the dynamical system (\ref{ODE4}) for $V$ and the PDE (\ref{W-PDE}) for $W$. 
The analysis in \cite{MM3} shows  the assumptions that $U\in W^{1,2}(\Omega)$ and that (\ref{DynSys}) is uniformly stable together imply that $V$ satisfies
\begin{equation}\label{est:C1-V}
\sup_{0<r<1}(|V(r)|+r|V'(r)|)\leq C
\end{equation}
and  that $W$ satisfies
\begin{equation}\label{est:M1p-W}
M_{1,p}(W,r)\leq  C\,\om(r) \,r \ \ \hbox{for}\ 0<r<1. 
\end{equation}
(In both (\ref{est:C1-V}) and (\ref{est:M1p-W}) the constants $C$ depend upon the $W^{1,2}$-norm of $U$, but not on $r$.)
Since $p>2$, we can use  Sobolev embedding to conclude that $|W(x,y)|r^{-1}\leq C\om(r)$ for $0<r<1$, i.e.\ $|W(x,y)|r^{-1}\to 0$ as $r=|{\bf x}|\to 0$. This shows that $W$ is differentiable at $0$.

To estimate $U_0$, we use (\ref{eq:ODE1}) and the estimates that we have obtained on
$V_i$ and $\nabla W$ to conclude
\[
\begin{aligned}
\left| U_0(r)-U_0(0)\right|=\left|\int_0^r U_0'(\rho)\,d\rho\right|&\leq 
C\,\omega(r)\,\int_0^r\left(\rho|V'(\rho)|+|V(\rho)|+|\nabla W|\right)\,d\rho\\
&\leq C\,\omega(r)\,r.
\end{aligned}
\]
But this implies that $U_0$ is differentiable at $x=0$ and $U_0'(0)=0$.

We have now shown  the assumption that (\ref{DynSys}) is uniformly stable is sufficient to show that our weak solution $U\in W^{1,2}(\Omega)$ satisfies
\[
\begin{aligned}
|U(x,y)-U(0)| & \leq |U_0(r)-U_0(0)|+|V_1(r)|\cdot|x|+|V_2(r)|\cdot|y|+|W(x,y)|  \leq C\,r.
\end{aligned}
\]
But this shows that $U$ is Lipschitz continuous at $x=0$, completing the proof of Theorem 1. 

For Theorem 2, we add the assumption that every solution of 
 (\ref{DynSys}) is asymptotically constant. The dynamical systems analysis of \cite{MM3} applied to (\ref{ODE5}) then shows that $\phi(t)\to\phi_\infty$ and $\psi(t)\to 0$ as $t\to\infty$. However, we can use (\ref{def:phi,psi}) to express $\phi,\psi$ in terms of $V$ and $V_t$:
 \[
 \begin{pmatrix} \phi \\ \psi \end{pmatrix}
 = \begin{pmatrix} \frac{1}{2}V-\frac{1}{4}V_t \\ \frac{1}{4}V_t \end{pmatrix} + O(\e(t)).
 \]
 Hence the conclusion $\psi\to 0$ implies $V_t\to 0$ as $t\to\infty$, in other words
 \begin{equation}\label{eq:rV->0}
 \lim_{r\to 0} r V'(r) =0.
 \end{equation}
 But (\ref{eq:rV->0}) implies that $V_1(r)x+V_2(r)y$ is differentiable at $0$.  Since we have already shown that $U_0$ and $W(x,y)$ are differentiable at $0$, we obtain the conclusion of Theorem 2.


\end{document}